\documentclass[12pt,a4paper]{article}
\usepackage[utf8]{inputenc}
\usepackage[T1]{fontenc}
\usepackage{amsmath}
\usepackage{amsthm}
\usepackage{amsfonts}
\usepackage{amssymb}
\usepackage{graphicx}
\usepackage{authblk}
\usepackage{hyperref}
\usepackage{tikz}
\usepackage{wrapfig}
\usepackage{caption}
\usepackage{subcaption}
\usepackage[capposition=top]{floatrow}
\usepackage[document]{ragged2e}
\usepackage[left=2.00cm, right=2.00cm, top=2.00cm, bottom=2.00cm]{geometry}

\newtheorem{theorem}{Theorem}

\newtheorem{observation}{Observation}

\title{\textbf{CONNECTIVITY OF ALTERNATING SIGN TRIANGLES}}
\author{SON NGUYEN}
\affil{University of Minnesota - Twin Cities\\\href{mailto:nguy4309@umn.edu}{nguy4309@umn.edu}}
\date{}

\begin{document}
\setlength{\abovedisplayskip}{0pt}
\setlength{\belowdisplayskip}{0pt}
\setlength{\abovedisplayshortskip}{0pt}
\setlength{\belowdisplayshortskip}{0pt}
	
	\maketitle
	
	\begin{abstract}
	\justify
		Alternating sign triangles were introduced by Carroll and Speyer in relation to cube recurrence, by analogy to alternating sign matrices for octahedron recurrence. In this paper, we prove the connectivity of alternating sign triangles, which is analogous to the connectivity of alternating sign matrices.
		
	\end{abstract}
	
	\section{Introduction}
	
	\justify
	The \textit{octahedron recurrence} is given by the initial conditions $f_{i,j,k}=x_{i,j,k}$ for $k = -1,0$ and
	
	\begin{align*}
	    f_{i,j,k-1}f_{i,j,k+1} &= f_{i-1,j,k}f_{i+1,j,k}+ \lambda f_{i,j-1,k}f_{i,j+1,k}. & (k\geq0)
	\end{align*}
	
	\justify
	An \textit{alternating sign matrix} is a square matrix that satisfies:
	
	\begin{enumerate}
	    \item all entries are $-1,0,\text{ or } 1$,
	    \item every row and column has sum $1$,
	    \item in every row and column the non-zero entries alternate in sign.
	\end{enumerate}
	
	\justify
	In \cite{robbins}, Robbins and Rumsey found that the exponents of the $x_{i,j,k}$ in any monomial in $f_{i_0,j_0,k_0}$ formed an alternating sign matrix. In \cite{brualdi}, the following theorem about alternating sign matrices was proven:
	
	\begin{theorem}
	    For every alternating sign matrices $A_{1}$ and $A_{2}$, $A_{1}$ can be obtained from $A_{2}$ by sequentially adding or subtracting the following sub-matrix:\\
	    
	    \centering
	    $\begin{bmatrix}
	        1 & -1 \\
	        -1 & 1
	    \end{bmatrix}$
	    
	\end{theorem}
	
	\justify
	In \cite{propp}, James Propp introduced another recurrence called the \textit{cube recurrence}, which is given by $f_{i,j,k} = x_{i,j,k}$ for $i+j+k=-1,0,1$ and
	
	\begin{align*}
	    f_{i,j,k}f_{i-1,j-1,k-1}&=f_{i-1,j,k}f_{i,j-1,k-1}+f_{i,j-1,k}f_{i-1,j,k-1}+f_{i,j,k-1}f_{i-1,j,k-1}. & (i+j+k>1)
	\end{align*}
	
	\justify
	In \cite{carroll}, Carroll and Speyer introduced a combinatorial object that describe the cube recurrence called \textit{grove}, which will be defined rigorously in section 2. Carroll and Speyer proved that $f_{0,0,0}$ is a sum of Laurent monomials in the variables $x_{i,j,k}$, and they observed that, in each monomial, the exponents of $x_{i,j,k}$ form an alternating sign triangles. Some properties of alternating sign triangles can be found in \cite{son}.
	
	In this paper, we will prove the following theorem, analogous to \textit{theorem 1}, for alternating sign triangles:
	
	\begin{theorem}
	    For every alternating sign triangles $A_{1}$ and $A_{2}$, $A_{1}$ can be obtained from $A_{2}$ by sequentially adding or subtracting one of the following three sub-triangles:
	\begin{figure}[h!]
		
     \centering
        \begin{tikzpicture}
            \draw (0,0) node[anchor=center] {-1};
            \draw (1,0) node[anchor=center] {1};
            \draw (0.5,-1) node[anchor=center] {0};
        \end{tikzpicture}
     \hspace{1cm}
        \begin{tikzpicture}
            \draw (0,0) node[anchor=center] {1};
            \draw (1,0) node[anchor=center] {0};
            \draw (0.5,-1) node[anchor=center] {-1};
        \end{tikzpicture}
     \hspace{1cm}
        \begin{tikzpicture}
            \draw (0,0) node[anchor=center] {0};
            \draw (1,0) node[anchor=center] {-1};
            \draw (0.5,-1) node[anchor=center] {1};
        \end{tikzpicture}
        
        \label{Sub-triangle}
        
    \end{figure}
	    
	\end{theorem}
    
    \justify
    \textbf{Acknowledgments.} I thank Pavlo Pylyavskyy for introducing me to this topic and for his helpful suggestions and continuous support.
	
	\section{Background}
	
	\justify
	In \cite{carroll}, Carroll and Speyer defined groves as follows.
	
	\justify
	Define the \textit{lower cone} of any $(i,j,k)\in \mathbb{Z}^{3}$ to be
	
	\begin{align*}
	    C(i,j,k) &= \left\{ (i',j',k')\in \mathbb{Z}^{3} | i'\leq i,j' \leq j,k'\leq k \right\}.
	\end{align*}
	
	\justify
    Let $\mathcal{L}\subseteq \mathbb{Z}^{3}$ be a subset such that, whenever $(i, j, k) \in \mathcal{L}, C(i, j, k) \subseteq \mathcal{L}$. Let $\mathcal{U} = \mathbb{Z}^{3}-\mathcal{L}$, and define the set of \textit{initial conditions}
    
    \begin{align*}
        \mathcal{I} = \left\{(i, j, k) \in \mathcal{L} | (i + 1, j + 1, k + 1) \in \mathcal{U}\right\}.
    \end{align*}
    
    \justify
    We also define a \textit{rhombus} to be any set of the form
    
    \begin{align*}
        r_{a}(i, j, k) &= \left\{(i, j, k),(i, j - 1, k),(i, j, k - 1),(i, j - 1, k - 1)\right\}\\
        r_{b}(i, j, k) &= \left\{(i, j, k),(i - 1, j, k),(i, j, k - 1),(i - 1, j, k - 1)\right\}\\
        r_{c}(i, j, k) &= \left\{(i, j, k),(i - 1, j, k),(i, j - 1, k),(i - 1, j - 1, k)\right\}
    \end{align*}
    
    \justify
    In addition, define the \textit{edges} of each rhombus to be the pairs
    
    \begin{align*}
        e_{a}(i, j, k) &= \left\{(i, j - 1, k),(i, j, k - 1)\right\} & e'_{a}(i, j, k) &= \left\{(i, j, k),(i, j - 1, k - 1)\right\}\\
        e_{b}(i, j, k) &= \left\{(i - 1, j, k),(i, j, k - 1)\right\} & e'_{b}(i, j, k) &= \left\{(i, j, k),(i - 1, j, k - 1)\right\}\\
        e_{c}(i, j, k) &= \left\{(i - 1, j, k),(i, j - 1, k)\right\} & e'_{c}(i, j, k) &= \left\{(i, j, k),(i - 1, j - 1, k)\right\}
    \end{align*}
    
    \justify
    Now suppose that $N$ is a cutoff for $\mathcal{I}$. Let $\mathcal{G}$ be a graph whose vertices are the points in $\mathcal{I}$ and edges are the edges of all rhombi occurring in $\mathcal{I}$. We define an $\mathcal{I}$-grove within radius $N$ to be a subgraph $G \subseteq \mathcal{G}$ with the following properties:
    
    \begin{itemize}
        \item (Completeness) the vertex set of $G$ is all of $\mathcal{I}$;
        \item (Complementarity) for every rhombus, exactly one of its two edges occurs in $G$;
        \item (Compactness) for every rhombus all of whose vertices satisfy $i+j+k<-N$, the short edge occurs in $G$;
        \item (Connectivity) every component of $G$ contains exactly one of the following sets of vertices, and conversely, each such set is contained in some component:
        \begin{itemize}
            \item $\left\{(0,p,q),(p,0,q)\right\},\left\{(p,q,0),(0,q,p)\right\}$, and $\left\{(q,0,p),(q,p,0)\right\}$ for all $p, q$ with $0>p>q$ and $p+q\in \{-N-1,-N-2\}$;
            \item $\left\{(0, p, p),(p, 0, p),(p, p, 0)\right\}$ for $2p\in \{-N-1,-N-2\}$;
            \item $\left\{(0, 0, q)\right\}, \left\{(0, q, 0)\right\}$, and $\left\{(q, 0, 0)\right\}$ for $q\leq-N-1$.
        \end{itemize}
    \end{itemize}
    
    \begin{figure}[h!]
        \centering
        \includegraphics[width=10cm]{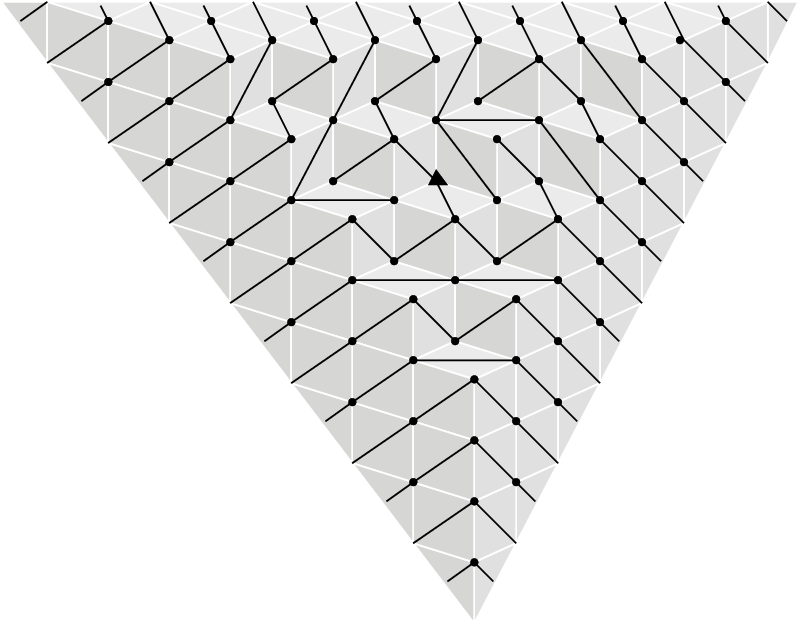}
        \caption{Example of a grove}
        \label{Figure 0}
        \floatfoot{Source: \cite{carroll}}
    \end{figure}
    
    \justify
    Carroll and Speyer also proved a bijection between groves and simplified groves. A simplified grove within radius $N$, where $N$ is a cutoff for $\mathcal{I}$ and furthermore is odd, is a subgraph $G'$ of $\mathcal{G}$ satisfying:
    
    \begin{itemize}
        \item (Vertex set) the vertex set of $G'$ is $\{(i,j,k)\in\mathcal{I} \text{ | } i+j+k\equiv0 \text{ mod } 2; i+j+k\geq-N-1\}$;
        \item (Acyclicity) $G'$ is acyclic;
        \item (Connectivity) every component of $G'$ contains exactly one of the following sets of vertices, and conversely, each such set is contained in some component:
        \begin{itemize}
            \item $\left\{(0,p,q),(p,0,q)\right\},\left\{(p,q,0),(0,q,p)\right\}$, and $\left\{(q,0,p),(q,p,0)\right\}$ for $p, q$ with $0>p>q$ and $p+q=-N-1$;
            \item $\left\{(0, \frac{-N-1}{2}, \frac{-N-1}{2}),(\frac{-N-1}{2}, 0, \frac{-N-1}{2}),(\frac{-N-1}{2}, \frac{-N-1}{2}, 0)\right\}$;
            \item $\left\{(0, 0, -N-1)\right\}, \left\{(0, -N-1, 0)\right\}$, and $\left\{(-N-1, 0, 0)\right\}$.
        \end{itemize}
    \end{itemize}
    
    \begin{figure}[h!]
        \centering
        \includegraphics[width=15cm]{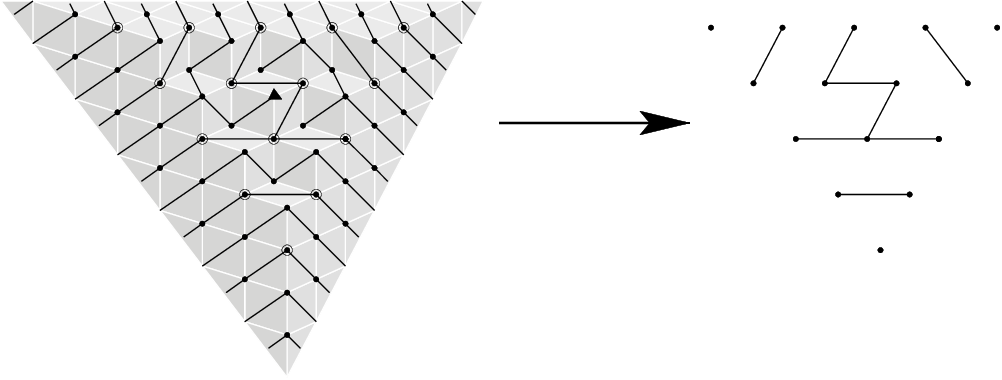}
        \caption{Example of a simplified grove}
        \label{Figure 1}
        \floatfoot{Source: \cite{carroll}}
    \end{figure}

	\justify
	\section{Definition}
	
	For our convenience, we will redefine a \textit{simplified grove of size $n$} to be a graph $G$ satisfying:
	
	\begin{itemize}
		\item (Vertex set) the vertex set of $G$ is $\{{(i,j)\in\mathbb{Z}^2\mid \lvert i+j\rvert\le n,j\leq 0,i+j\equiv n\bmod{2}}\}$;
		\item (Acyclicity) $G$ is acyclic;
		\item (Connectivity) the boundary vertices can be partitioned into the following sets so that each component of $G$ contains exactly one set, and conversely, each set is contained in some component:
		\begin{itemize}
		    \item $\{(-n,0)\},\{(n,0)\}$ and $\{(0,-n)\}$,
		    \item (west pairs) $\{(-i,0),(\frac{-n-i}{2},\frac{-n+i}{2})\}$ for $0<i<n, i\equiv n\bmod{2}$,
		    \item (east pairs) $\{(i,0),(\frac{n+i}{2},\frac{-n+i}{2})\}$ for $0<i<n, i\equiv n\bmod{2}$,
		    \item (south pairs) $\{(-n+i,-i),(n-i,-i\}$ for $\frac{n}{2}<i<n$,
		    \item (middle triplet) $\{(0,0),(-\frac{n}{2},-\frac{n}{2}),(\frac{n}{2},-\frac{n}{2})\}$ if $n$ is even.
		\end{itemize}
	\end{itemize}

	It can be checked that this new definition gives the same set of simplified groves as Carroll's. We also define a \textit{upward triangle of size $i$} to be a triangle whose vertices are $(a,b),(a-i,b-i),(a+i,b-i)$. Similarly, define an \textit{downward triangle of size $i$} to be a triangle whose vertices are $(a,b),(a-i,b+i),(a+i,b+i)$. For simplicity, we will refer to downward and upward triangle of size $1$ as downward and upward triangle respectively.
	
	Now assign to each downward triangle a number $1-e$ where $e$ is the number of edges in that triangle. From the acyclicity condition, we have the number of edge in every triangle is less than $3$; hence, the number assigned to each downward triangle is either $-1,0$ or $1$. Define an \textit{alternating sign triangle} to be a configuration of numbers generated by our process. It can be checked that this definition gives the same set of alternating sign triangles as Carroll's proposed definition does.
    
    \begin{figure}[h!]
        \centering
    
    \begin{tikzpicture}
        \filldraw[black] (0,0) circle (2pt);
        \filldraw[black] (1,0) circle (2pt);
        \filldraw[black] (2,0) circle (2pt);
        \filldraw[black] (3,0) circle (2pt);
        \filldraw[black] (4,0) circle (2pt);
        \filldraw[black] (0.5,-1) circle (2pt);
        \filldraw[black] (1.5,-1) circle (2pt);
        \filldraw[black] (2.5,-1) circle (2pt);
        \filldraw[black] (3.5,-1) circle (2pt);
        \filldraw[black] (1,-2) circle (2pt);
        \filldraw[black] (2,-2) circle (2pt);
        \filldraw[black] (3,-2) circle (2pt);
        \filldraw[black] (1.5,-3) circle (2pt);
        \filldraw[black] (2.5,-3) circle (2pt);
        \filldraw[black] (2,-4) circle (2pt);
        
        \draw[black, thick] (1,0) -- (0.5,-1);
        \draw[black, thick] (2,0) -- (2.5,-1);
        \draw[black, thick] (2.5,-1) -- (1.5,-1);
        \draw[black, thick] (1.5,-1) -- (1,-2);
        \draw[black, thick] (1.5,-1) -- (2,-2);
        \draw[black, thick] (2,-2) -- (3,-2);
        \draw[black, thick] (3,0) -- (3.5,-1);
        \draw[black, thick] (1.5,-3) -- (2.5,-3);
        
        \draw [-latex, thick](6,-2) -- (8,-2);
        
        \draw (10,-0.5) node[anchor=center] {0};
        \draw (11,-0.5) node[anchor=center] {1};
        \draw (12,-0.5) node[anchor=center] {0};
        \draw (13,-0.5) node[anchor=center] {0};
        \draw (10.5,-1.5) node[anchor=center] {0};
        \draw (11.5,-1.5) node[anchor=center] {-1};
        \draw (12.5,-1.5) node[anchor=center] {1};
        \draw (11,-2.5) node[anchor=center] {1};
        \draw (12,-2.5) node[anchor=center] {0};
        \draw (11.5,-3.5) node[anchor=center] {0};
        
    \end{tikzpicture}
        \caption{Getting alternating sign triangle of size 4 from a grove of size 4}
        \label{Figure 2}
    \end{figure}
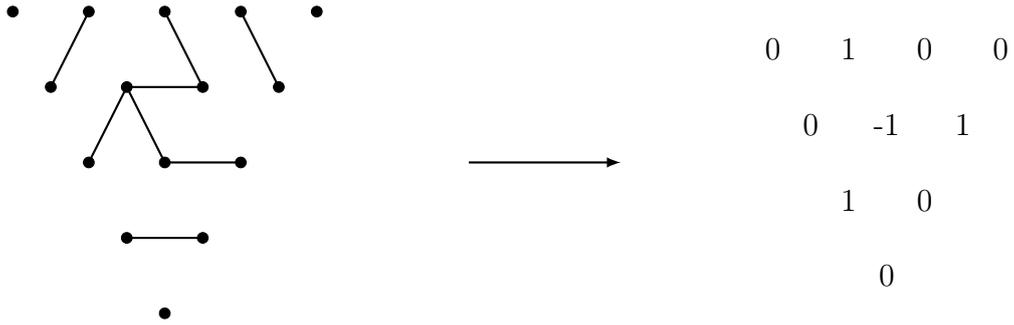
	
	\newpage
	\justify
    Now we will introduce some definitions that will be useful for our proof. First of all, define the \textit{target grove of size n} to be the following grove:
    
    \begin{figure}[h!]
        \centering
    
    \begin{tikzpicture}
        \filldraw[black] (0,0) circle (2pt);
        \filldraw[black] (1,0) circle (2pt);
        \filldraw[black] (2,0) circle (2pt);
        \filldraw[black] (3,0) circle (2pt);
        \filldraw[black] (4,0) circle (2pt);
        \filldraw[black] (5,0) circle (2pt);
        \filldraw[black] (0.5,-1) circle (2pt);
        \filldraw[black] (1.5,-1) circle (2pt);
        \filldraw[black] (2.5,-1) circle (2pt);
        \filldraw[black] (3.5,-1) circle (2pt);
        \filldraw[black] (4.5,-1) circle (2pt);
        \filldraw[black] (1,-2) circle (2pt);
        \filldraw[black] (2,-2) circle (2pt);
        \filldraw[black] (3,-2) circle (2pt);
        \filldraw[black] (4,-2) circle (2pt);
        \filldraw[black] (1.5,-3) circle (2pt);
        \filldraw[black] (2.5,-3) circle (2pt);
        \filldraw[black] (3.5,-3) circle (2pt);
        \filldraw[black] (2,-4) circle (2pt);
        \filldraw[black] (3,-4) circle (2pt);
        \filldraw[black] (2.5,-5) circle (2pt);
        
        \draw[black, thick, dotted] (2.8,0) -- (2.2,0);
        \draw[black, thick, dotted] (3.3,-1) -- (2.7,-1);
        \draw[black, thick, dotted] (3.8,-2) -- (3.2,-2);
        \draw[black, thick, dotted] (4,-2) -- (3.5,-3);
        \draw[black, thick, dotted] (3,-2) -- (3.5,-3);
        \draw[black, thick, dotted] (2,-2) -- (2.5,-3);
        \draw[black, thick, dotted] (1,-2) -- (1.5,-3);
        
        \draw[black, thick] (1,0) -- (1.5,-1);
        \draw[black, thick] (1.5,-1) -- (0.5,-1);
        \draw[black, thick] (2,0) -- (3,-2);
        \draw[black, thick] (3,-2) -- (1,-2);
        \draw[black, thick] (3,0) -- (4,-2);
        \draw[black, thick] (4,0) -- (4.5,-1);
        \draw[black, thick] (3.5,-3) -- (1.5,-3);
        \draw[black, thick] (3,-4) -- (2,-4);
        
        \draw (-0.5,-2) node[anchor=center] {$W$};
        \draw (5.4,-2) node[anchor=center] {$E$};
        \draw (0.1,-3) node[anchor=center] {$S$};
        
        \draw [-latex, thick](0,-2) -- (0.8,-2);
        \draw [-latex, thick](5,-2) -- (4.2,-2);
        \draw [-latex, thick](0.5,-3) -- (1.3,-3);
        
    \end{tikzpicture}
        \caption{Target grove}
        \label{Figure 3}
    \end{figure}
    
    \justify
    We also say two vertices or two edges to be in the same group if they are in the same component in the target grove. Additionally, we define group $W,S,E$ as in Figure \ref{Figure 3} as they play an important part in our proof. Next, for every grove $G$, define its \textit{difference grove} to be the graph in which:
    
    \begin{itemize}
        \item Edges that appear in $G$ but not in the target grove are colored red.
        \item Edges that appear the target grove but not in $G$ are colored black.
        \item Edges that appear in both $G$ and the target grove are colored blue.
    \end{itemize}
    
    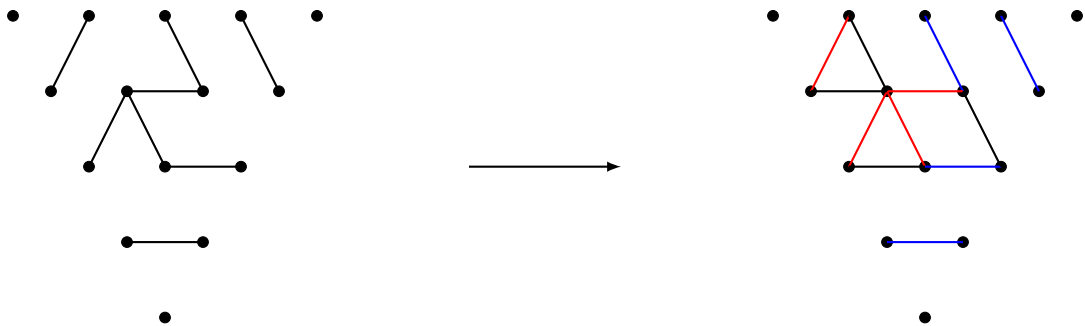
\begin{figure}[h!]
        \centering
    
    \begin{tikzpicture}
        \filldraw[black] (0,0) circle (2pt);
        \filldraw[black] (1,0) circle (2pt);
        \filldraw[black] (2,0) circle (2pt);
        \filldraw[black] (3,0) circle (2pt);
        \filldraw[black] (4,0) circle (2pt);
        \filldraw[black] (0.5,-1) circle (2pt);
        \filldraw[black] (1.5,-1) circle (2pt);
        \filldraw[black] (2.5,-1) circle (2pt);
        \filldraw[black] (3.5,-1) circle (2pt);
        \filldraw[black] (1,-2) circle (2pt);
        \filldraw[black] (2,-2) circle (2pt);
        \filldraw[black] (3,-2) circle (2pt);
        \filldraw[black] (1.5,-3) circle (2pt);
        \filldraw[black] (2.5,-3) circle (2pt);
        \filldraw[black] (2,-4) circle (2pt);
        
        \draw[black, thick] (1,0) -- (0.5,-1);
        \draw[black, thick] (2,0) -- (2.5,-1);
        \draw[black, thick] (2.5,-1) -- (1.5,-1);
        \draw[black, thick] (1.5,-1) -- (1,-2);
        \draw[black, thick] (1.5,-1) -- (2,-2);
        \draw[black, thick] (2,-2) -- (3,-2);
        \draw[black, thick] (3,0) -- (3.5,-1);
        \draw[black, thick] (1.5,-3) -- (2.5,-3);
        
        \draw [-latex, thick](6,-2) -- (8,-2);
        
        \filldraw[black] (10,0) circle (2pt);
        \filldraw[black] (11,0) circle (2pt);
        \filldraw[black] (12,0) circle (2pt);
        \filldraw[black] (13,0) circle (2pt);
        \filldraw[black] (14,0) circle (2pt);
        \filldraw[black] (10.5,-1) circle (2pt);
        \filldraw[black] (11.5,-1) circle (2pt);
        \filldraw[black] (12.5,-1) circle (2pt);
        \filldraw[black] (13.5,-1) circle (2pt);
        \filldraw[black] (11,-2) circle (2pt);
        \filldraw[black] (12,-2) circle (2pt);
        \filldraw[black] (13,-2) circle (2pt);
        \filldraw[black] (11.5,-3) circle (2pt);
        \filldraw[black] (12.5,-3) circle (2pt);
        \filldraw[black] (12,-4) circle (2pt);

        \draw[red, thick] (11,0) -- (10.5,-1);
        \draw[blue, thick] (12,0) -- (12.5,-1);
        \draw[red, thick] (12.5,-1) -- (11.5,-1);
        \draw[red, thick] (11.5,-1) -- (11,-2);
        \draw[red, thick] (11.5,-1) -- (12,-2);
        \draw[blue, thick] (12,-2) -- (13,-2);
        \draw[blue, thick] (13,0) -- (13.5,-1);
        \draw[blue, thick] (11.5,-3) -- (12.5,-3);
        \draw[black, thick] (12.5,-1) -- (13,-2);
        \draw[black, thick] (11,-2) -- (12,-2);
        \draw[black, thick] (11,0) -- (11.5,-1);
        \draw[black, thick] (10.5,-1) -- (11.5,-1);
        
    \end{tikzpicture}
        \caption{A grove and its difference grove}
        \label{Figure 4}
    \end{figure}
    
    \newpage
    \justify
    Lastly, we define a \textit{spin} to be one of the following changes:\\
    
    \centering
    \begin{tikzpicture}
    
        \filldraw[black] (0.5,0) circle (2pt);
        \filldraw[black] (1.5,0) circle (2pt);
        \filldraw[black] (0,-1) circle (2pt);
        \filldraw[black] (1,-1) circle (2pt);
        \filldraw[black] (2,-1) circle (2pt);
        \filldraw[black] (0.5,-2) circle (2pt);
        \filldraw[black] (1.5,-2) circle (2pt);
        
        \draw[black, thick] (1,-1) -- (0.5,0);

        \draw [-latex, thick](3,-0.8) -- (4,-0.8);
        \draw [latex-, thick](3,-1.2) -- (4,-1.2);
    
        \filldraw[black] (5.5,0) circle (2pt);
        \filldraw[black] (6.5,0) circle (2pt);
        \filldraw[black] (5,-1) circle (2pt);
        \filldraw[black] (6,-1) circle (2pt);
        \filldraw[black] (7,-1) circle (2pt);
        \filldraw[black] (5.5,-2) circle (2pt);
        \filldraw[black] (6.5,-2) circle (2pt);
        
        \draw[black, thick] (6,-1) -- (6.5,0);

        \draw [-latex, thick](8,-0.8) -- (9,-0.8);
        \draw [latex-, thick](8,-1.2) -- (9,-1.2);
    
        \filldraw[black] (10.5,0) circle (2pt);
        \filldraw[black] (11.5,0) circle (2pt);
        \filldraw[black] (10,-1) circle (2pt);
        \filldraw[black] (11,-1) circle (2pt);
        \filldraw[black] (12,-1) circle (2pt);
        \filldraw[black] (10.5,-2) circle (2pt);
        \filldraw[black] (11.5,-2) circle (2pt);
        
        \draw[black, thick] (11,-1) -- (12,-1);

        \draw [-latex, thick](10.8,-3) -- (10.8,-4);
        \draw [latex-, thick](11.2,-3) -- (11.2,-4);
    
        \filldraw[black] (0.5,-5) circle (2pt);
        \filldraw[black] (1.5,-5) circle (2pt);
        \filldraw[black] (0,-6) circle (2pt);
        \filldraw[black] (1,-6) circle (2pt);
        \filldraw[black] (2,-6) circle (2pt);
        \filldraw[black] (0.5,-7) circle (2pt);
        \filldraw[black] (1.5,-7) circle (2pt);
        
        \draw[black, thick] (1,-6) -- (0,-6);

        \draw [-latex, thick](3,-5.8) -- (4,-5.8);
        \draw [latex-, thick](3,-6.2) -- (4,-6.2);
    
        \filldraw[black] (5.5,-5) circle (2pt);
        \filldraw[black] (6.5,-5) circle (2pt);
        \filldraw[black] (5,-6) circle (2pt);
        \filldraw[black] (6,-6) circle (2pt);
        \filldraw[black] (7,-6) circle (2pt);
        \filldraw[black] (5.5,-7) circle (2pt);
        \filldraw[black] (6.5,-7) circle (2pt);
        
        \draw[black, thick] (6,-6) -- (5.5,-7);

        \draw [-latex, thick](8,-5.8) -- (9,-5.8);
        \draw [latex-, thick](8,-6.2) -- (9,-6.2);
    
        \filldraw[black] (10.5,-5) circle (2pt);
        \filldraw[black] (11.5,-5) circle (2pt);
        \filldraw[black] (10,-6) circle (2pt);
        \filldraw[black] (11,-6) circle (2pt);
        \filldraw[black] (12,-6) circle (2pt);
        \filldraw[black] (10.5,-7) circle (2pt);
        \filldraw[black] (11.5,-7) circle (2pt);
        
        \draw[black, thick] (11,-6) -- (11.5,-7);

        \draw [-latex, thick](0.8,-3) -- (0.8,-4);
        \draw [latex-, thick](1.2,-3) -- (1.2,-4);
        
        \draw (3.5,-0.5) node[anchor=center] {(1)};
        \draw (8.5,-0.5) node[anchor=center] {(2)};
        \draw (11.7,-3.5) node[anchor=center] {(3)};
        \draw (8.5,-5.5) node[anchor=center] {(4)};
        \draw (3.5,-5.5) node[anchor=center] {(5)};
        \draw (0.3,-3.5) node[anchor=center] {(6)};
        
    \end{tikzpicture}
    
    \justify
    We also call the middle vertex above the \textit{pivot} of the spin. It can be seen that spins of type $1,3,5$ does not change the corresponding alternating sign triangles while spins of type $2,4,6$ is equivalent to adding or subtracting one of the three sub-triangles in Theorem 2 to the corresponding alternating sign triangles. It can also be seen that if a grove $G'$ can be obtained from another grove $G$ through a sequence of spins, then $G$ can also be obtained from $G'$ through a sequence of spins. Hence, instead of proving Theorem 2, we will prove the following theorem:
    
    \begin{theorem}
        For every grove $G$, we can obtain the target grove after a sequence of spins.
    \end{theorem}
    
    \section{Observations}
    \justify
    Before proving Theorem 3, we will make some trivial observations:
    
    \begin{observation}
        Between any two consecutive black edges in a group, and between any black edge and the boundary vertices, there is at least one red edge.
    \end{observation}
    
    \begin{observation}
        If an interior vertex has degree 1, we can spin freely with that vertex as the pivot.
    \end{observation}
    
    \begin{observation}
        We can move any black edge to its nearest red edge.
    \end{observation}
    
    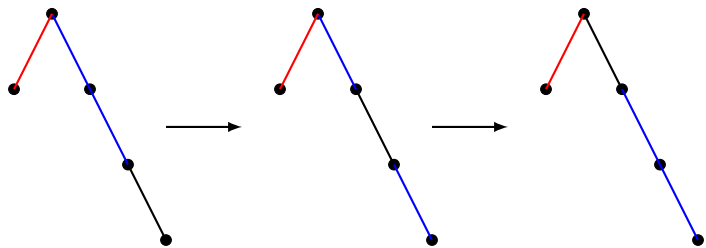
\begin{figure}[h!]
        \centering
        \begin{tikzpicture}
        \filldraw[black] (0,0) circle (2pt);
        \filldraw[black] (0.5,-1) circle (2pt);
        \filldraw[black] (-0.5,-1) circle (2pt);
        \filldraw[black] (1,-2) circle (2pt);
        \filldraw[black] (1.5,-3) circle (2pt);
        
        \draw[black, thick] (1.5,-3) -- (1,-2);
        \draw[blue, thick] (0,0) -- (1,-2);
        \draw[red, thick] (0,0) -- (-0.5,-1);
        
        \draw [-latex, thick](1.5,-1.5) -- (2.5,-1.5);
        
        \filldraw[black] (3.5,0) circle (2pt);
        \filldraw[black] (4,-1) circle (2pt);
        \filldraw[black] (3,-1) circle (2pt);
        \filldraw[black] (4.5,-2) circle (2pt);
        \filldraw[black] (5,-3) circle (2pt);
        
        \draw[black, thick] (4,-1) -- (4.5,-2);
        \draw[blue, thick] (4,-1) -- (3.5,0);
        \draw[blue, thick] (4.5,-2) -- (5,-3);
        \draw[red, thick] (3.5,0) -- (3,-1);
        
        \draw [-latex, thick](5,-1.5) -- (6,-1.5);
        
        \filldraw[black] (7,0) circle (2pt);
        \filldraw[black] (7.5,-1) circle (2pt);
        \filldraw[black] (6.5,-1) circle (2pt);
        \filldraw[black] (8,-2) circle (2pt);
        \filldraw[black] (8.5,-3) circle (2pt);
        
        \draw[black, thick] (7.5,-1) -- (7,0);
        \draw[blue, thick] (7.5,-1) -- (8.5,-3);
        \draw[red, thick] (7,0) -- (6.5,-1);
        
        \end{tikzpicture}
        \caption{Observation 3}
        \label{Figure 5}
    \end{figure}
    
    \section{Proof}
    
    \justify
    We will prove that if there is some black edges in the difference grove, we can reduce the number of black edges.
    
    \textit{Case 1:} If there is no group with more than one black edge, then, we choose a group with one black edge. Clearly, the path connecting the two boundary vertices of this group must contain a segment of red edges. Moreover, the black edge lies between the two endpoints of the segment. Because there is only one black edge, every vertices in the region bounded by the red segment and the group edges must be in the same component with the group's vertices. This means that we can move the black edge to one of the endpoints of the red segment and then, with the endpoint as the pivot, spin the red edge to the place of the black edge, and the number of black edge is reduced.
    
    \begin{figure}[h!]
        \centering
        \begin{tikzpicture}
        \filldraw[black] (0,0) circle (2pt);
        \filldraw[black] (0.5,-1) circle (2pt);
        \filldraw[black] (1,-2) circle (2pt);
        \filldraw[black] (1.5,-3) circle (2pt);
        \filldraw[black] (0.5,-3) circle (2pt);
        \filldraw[black] (-0.5,-3) circle (2pt);
        \filldraw[black] (-0.5,-1) circle (2pt);
        \filldraw[black] (-1,-2) circle (2pt);
        
        \draw[black, thick] (1.5,-3) -- (1,-2);
        \draw[red, thick] (0,0) -- (-1,-2);
        \draw[red, thick] (-0.5,-3) -- (-1,-2);
        \draw[blue, thick] (1.5,-3) -- (-0.5,-3);
        \draw[blue, thick] (1,-2) -- (0,0);
        
        \draw [-latex, thick](1.5,-1.5) -- (2.5,-1.5);
        
        \filldraw[black] (4,0) circle (2pt);
        \filldraw[black] (4.5,-1) circle (2pt);
        \filldraw[black] (5,-2) circle (2pt);
        \filldraw[black] (5.5,-3) circle (2pt);
        \filldraw[black] (4.5,-3) circle (2pt);
        \filldraw[black] (3.5,-3) circle (2pt);
        \filldraw[black] (3.5,-1) circle (2pt);
        \filldraw[black] (3,-2) circle (2pt);
        
        \draw[black, thick] (3.5,-3) -- (4.5,-3);
        \draw[red, thick] (4,0) -- (3,-2);
        \draw[red, thick] (3.5,-3) -- (3,-2);
        \draw[blue, thick] (5.5,-3) -- (4.5,-3);
        \draw[blue, thick] (5.5,-3) -- (4,0);
        
        \draw [-latex, thick](5.5,-1.5) -- (6.5,-1.5);
        
        \filldraw[black] (8,0) circle (2pt);
        \filldraw[black] (8.5,-1) circle (2pt);
        \filldraw[black] (9,-2) circle (2pt);
        \filldraw[black] (9.5,-3) circle (2pt);
        \filldraw[black] (8.5,-3) circle (2pt);
        \filldraw[black] (7.5,-3) circle (2pt);
        \filldraw[black] (7.5,-1) circle (2pt);
        \filldraw[black] (7,-2) circle (2pt);
        
        \draw[red, thick] (8,0) -- (7,-2);
        \draw[blue, thick] (9.5,-3) -- (7.5,-3);
        \draw[blue, thick] (9.5,-3) -- (8,0);
        
        \end{tikzpicture}
        \caption{Case 1}
        \label{Figure 6}
    \end{figure}
    
    \justify
    \textit{Case 2:} If there are some groups with more than one black edge, then, we choose one of such groups and consider a pair of consecutive black edges. By \textit{observation 1}, there have to be at least one red edge between the two black edges. If there is only one red edge, by \textit{observation 3}, we can move the two black edges to the red edge. Now the shared vertex of the two black edges and the red edge has degree 1; therefore, by \textit{observation 2}, we can spin with that vertex as the pivot to move the red edge to the place of one of the black edges, and the number of black edges is reduced.
    
    On the other hand, if there are more than one red edge, we consider two adjacent red edges and consider two sub cases:
    
    \textit{Case 2a:} If the two chosen red edges are connected to the same component, then in the other component, between the two red edges, there must be some black edges since otherwise we will have a cycle. If there is only one black edge, then similar to \textit{case 1}, we can reduce the number of black edges.\\
    
    \centering
    \begin{tikzpicture}
        \filldraw[black] (0,0) circle (2pt);
        \filldraw[black] (0.5,-1) circle (2pt);
        \filldraw[black] (1,-2) circle (2pt);
        \filldraw[black] (1.5,-3) circle (2pt);
        \filldraw[black] (0.5,-3) circle (2pt);
        \filldraw[black] (0,-2) circle (2pt);
        \filldraw[black] (-0.5,-1) circle (2pt);
        \filldraw[black] (-1,0) circle (2pt);
        
        \draw[blue, thick] (0,0) -- (1.5,-3);
        \draw[red, thick] (0,0) -- (-1,0);
        \draw[red, thick] (0.5,-3) -- (1.5,-3);
        \draw[blue, thick] (-0.5,-1) -- (-1,0);
        \draw[blue, thick] (0.5,-3) -- (0,-2);
        \draw[black, thick] (-0.5,-1) -- (0,-2);
        
        \draw [-latex, thick](1.75,-1.5) -- (2.75,-1.5);
        
        \filldraw[black] (4,0) circle (2pt);
        \filldraw[black] (4.5,-1) circle (2pt);
        \filldraw[black] (5,-2) circle (2pt);
        \filldraw[black] (5.5,-3) circle (2pt);
        \filldraw[black] (4.5,-3) circle (2pt);
        \filldraw[black] (4,-2) circle (2pt);
        \filldraw[black] (3.5,-1) circle (2pt);
        \filldraw[black] (3,0) circle (2pt);
        
        \draw[blue, thick] (4,0) -- (5.5,-3);
        \draw[red, thick] (4,0) -- (3,0);
        \draw[red, thick] (4.5,-3) -- (5.5,-3);
        \draw[blue, thick] (4,-2) -- (3,0);
        \draw[black, thick] (4.5,-3) -- (4,-2);
        
        \draw [-latex, thick](5.75,-1.5) -- (6.75,-1.5);
        
        \filldraw[black] (8,0) circle (2pt);
        \filldraw[black] (8.5,-1) circle (2pt);
        \filldraw[black] (9,-2) circle (2pt);
        \filldraw[black] (9.5,-3) circle (2pt);
        \filldraw[black] (8.5,-3) circle (2pt);
        \filldraw[black] (8,-2) circle (2pt);
        \filldraw[black] (7.5,-1) circle (2pt);
        \filldraw[black] (7,0) circle (2pt);
        
        \draw[blue, thick] (8,0) -- (9.5,-3);
        \draw[red, thick] (8,0) -- (7,0);
        \draw[blue, thick] (8.5,-3) -- (7,0);
        
    \end{tikzpicture}
    
    \justify
    If there are more than one black edge, consider two consecutive ones, then by \textit{observation 1}, there have to be some red edges between them. Because of our choice of red edges at the beginning of this case, these red edges cannot be connected with our original group; hence, they have to be connected to other groups. If there is only one red edge, then similar as above, we can reduce the number of red edges. If there are more than one red edge, and they are connected to the same group, then we get back to this case. If they are connected to different groups, then we get the following bad scenario, which will be considered last.\\
    
    \centering
    \begin{tikzpicture}
        \filldraw[black] (0,0) circle (2pt);
        \filldraw[black] (0.5,-1) circle (2pt);
        \filldraw[black] (1,-2) circle (2pt);
        \filldraw[black] (1.5,-3) circle (2pt);

        \filldraw[black] (0.5,-3) circle (2pt);
        \filldraw[black] (0,-2) circle (2pt);
        \filldraw[black] (-0.5,-1) circle (2pt);
        \filldraw[black] (-1,0) circle (2pt);
        \filldraw[black] (-0.5,-3) circle (2pt);
        \filldraw[black] (-1.5,-3) circle (2pt);
        \filldraw[black] (-2.5,-3) circle (2pt);
        
        \filldraw[black] (0,-4) circle (2pt);
        \filldraw[black] (-1,-4) circle (2pt);
        \filldraw[black] (-2,-4) circle (2pt);
        \filldraw[black] (1,-4) circle (2pt);
        
        \draw[blue, thick] (0.5,-3) -- (-2.5,-3);
        \draw[blue, thick] (0.5,-3) -- (-1,0);
        \draw[blue, thick] (1.5,-3) -- (0,0);
        \draw[blue, thick] (-2,-4) -- (1,-4);
        \draw[red, thick] (-0.5,-3) -- (-1,-4);
        \draw[red, thick] (0,-2) -- (1,-2);
        
    \end{tikzpicture}
    
    \justify
    \textit{Case 2b:} If the two chosen red edges are connected to different components, then in at least one of the two other components, there have to be black edges on both sides of the red edge because otherwise the boundaries vertices of the two components will be connected. Then if between the two black edges there are only one red edge, then similar as above, we can reduce to number of black edges. If there are more than one red edge, then we proceed similarly as above. Consider two adjacent red edges, if they are connected to the same component, then we are back to \textit{case 2a}, else, we are back to this case.
    
    Now let us take a moment to review what we have got. We now know that if we get \textit{case 2b}, we will either finish, or get back to \textit{case 2b} again, or get to \textit{case 2a}. If we get \textit{case 2a}, we will either finish, or get back to \textit{case 2a} again, or get to the \textit{bad scenario}. However, this process cannot go on forever because we only have a finite number of red edges, so the process will eventually stops. Therefore, eventually, we will either finish, or get to the \textit{bad scenario}. Hence, if we can solve for the \textit{bad scenario}, the proof is complete.
    
    Now we will show that the bad scenario is actually not that bad. To solve for the bad scenario, we need to see when the bad scenario arises. It arises when we have three groups adjacent to each other, and the only place where that happens is at group $W,E,S$ that we defined earlier. Additionally, it only arises when, in our process, we go from group $W$ to group $E$ and $S$ in \textit{case 2a}.\\
    
    \centering
    \begin{tikzpicture}
        \filldraw[black] (0,0) circle (2pt);
        \filldraw[black] (0.5,-1) circle (2pt);
        \filldraw[black] (1,-2) circle (2pt);
        \filldraw[black] (1.5,-3) circle (2pt);

        \filldraw[black] (0.5,-3) circle (2pt);
        \filldraw[black] (0,-2) circle (2pt);
        \filldraw[black] (-0.5,-1) circle (2pt);
        \filldraw[black] (-1,0) circle (2pt);
        \filldraw[black] (-0.5,-3) circle (2pt);
        \filldraw[black] (-1.5,-3) circle (2pt);
        \filldraw[black] (-2.5,-3) circle (2pt);
        
        \filldraw[black] (0,-4) circle (2pt);
        \filldraw[black] (-1,-4) circle (2pt);
        \filldraw[black] (-2,-4) circle (2pt);
        \filldraw[black] (1,-4) circle (2pt);
        
        \draw[blue, thick] (0.5,-3) -- (-2.5,-3);
        \draw[blue, thick] (0.5,-3) -- (-1,0);
        \draw[blue, thick] (1.5,-3) -- (0,0);
        \draw[blue, thick] (-2,-4) -- (1,-4);
        \draw[red, thick] (-0.5,-3) -- (-1,-4);
        \draw[red, thick] (0,-2) -- (1,-2);
        
        \draw (-1.5,-2) node[anchor=center] {$W$};
        \draw (2.4,-2) node[anchor=center] {$E$};
        \draw (-3.4,-4) node[anchor=center] {$S$};
        
        \draw [-latex, thick](-1,-2) -- (-0.2,-2);
        \draw [-latex, thick](2,-2) -- (1.2,-2);
        \draw [-latex, thick](-3,-4) -- (-2.2,-4);
        
    \end{tikzpicture}
    
    \justify
    Now if we look closely, we can see that we can apply the same reasoning as in \textit{case 2b} for group $S$ and $E$, that is, on at least one of the two groups, there have to be black edges on both sides of the red edge. Then we can proceed as we did with \textit{case 2b}, and we will have that we either finish, or get to \textit{case 2a}, or get back to \textit{case 2b}, and the process goes on as above. However, we can see that if we continue getting \textit{case 2b}, our process will only go either eastward or southward. If at any point we get \textit{case 2a}, our process will only go westward until it stops. That means that we will never get to the bad scenario again, which means that when our process stops, and we know that it will stop eventually, we will finish. This completes our proof.
    
    \section{Discussion}
    
    \justify
    Now looking back at our proof, we can see that we use only two types of spins: spins with the pivot having degree one, and spins among vertices in the same component as in \textit{case 1}. In the first type, which direction we spin does not matter, so we can always spin in, say, clockwise direction. In the second type, notice that one of the endpoints of the red segment is to the left of the black edge, so if we use that endpoint as the pivot, we can spin in the clockwise direction only. Therefore, we have the following stronger theorem of \textit{theorem 3}:
    
    \begin{theorem}
        For every grove $G$, we can obtain the target grove after a sequence of spins in the clockwise direction.
    \end{theorem}

    \justify
    Now let us define the \textit{identity triangle}, analogous to \textit{identity matrix}, to be the following triangle:
    
	\begin{figure}[h!]
		
     \centering
        \begin{tikzpicture}
            \draw (0,0) node[anchor=center] {1};
            \draw (1,0) node[anchor=center] {0};
            \draw (2,0) node[anchor=center] {0};
            \draw (3,0) node[anchor=center] {0};
            \draw (4,0) node[anchor=center] {0};
            \draw (6,0) node[anchor=center] {0};
            \draw (0.5,-1) node[anchor=center] {0};
            \draw (1.5,-1) node[anchor=center] {1};
            \draw (2.5,-1) node[anchor=center] {0};
            \draw (3.5,-1) node[anchor=center] {0};
            \draw (5.5,-1) node[anchor=center] {0};
            \draw (1,-2) node[anchor=center] {0};
            \draw (2,-2) node[anchor=center] {0};
            \draw (3,-2) node[anchor=center] {1};
            \draw (5,-2) node[anchor=center] {0};
            \draw (3.5,-5) node[anchor=center] {0};
            \draw (2.5,-5) node[anchor=center] {0};
            \draw (3,-6) node[anchor=center] {0};
            
            \draw[black, dashed] (4.5,0) -- (5.5,0);
            \draw[black, dashed] (4,-1) -- (5,-1);
            \draw[black, dashed] (4.5,-2) -- (3.5,-2);
            \draw[black, dashed] (1.5,-3) -- (4.5,-3);
            \draw[black, dashed] (2,-4) -- (4,-4);
        \end{tikzpicture}
    \end{figure}
    
    Clearly, the target grove gives the identity triangle; therefore, from \textit{theorem 4}, we have the following theorem:
    
	\begin{theorem}
	    From every alternating sign triangle, we can obtain the identity triangle by sequentially adding one of the following three sub-triangles:
	    
	\begin{figure}[h!]
		
     \centering
        \begin{tikzpicture}
            \draw (0,0) node[anchor=center] {-1};
            \draw (1,0) node[anchor=center] {1};
            \draw (0.5,-1) node[anchor=center] {0};
        \end{tikzpicture}
     \hspace{1cm}
        \begin{tikzpicture}
            \draw (0,0) node[anchor=center] {1};
            \draw (1,0) node[anchor=center] {0};
            \draw (0.5,-1) node[anchor=center] {-1};
        \end{tikzpicture}
     \hspace{1cm}
        \begin{tikzpicture}
            \draw (0,0) node[anchor=center] {0};
            \draw (1,0) node[anchor=center] {-1};
            \draw (0.5,-1) node[anchor=center] {1};
        \end{tikzpicture}
        
        \label{Sub-triangle}
        
    \end{figure}
	    
	\end{theorem}
	
    \justify
    There are also bijections between \textit{alternating sign matrix} and some other interesting objects, one of which is the \textit{corner-sum matrix}, which was introduced in \cite{robbins}. Corner-sum matrices are matrices in which the last row and last column consist of the numbers from 1 to n (in order), and within each row and column, each entry is either equal to, or one more than, the preceding entry. By adding (or subtracting) the sub-matrix in \textit{theorem 1} to an alternating sign matrix, we are increasing (or decreasing) the top left entry of the sub-matrix in the corresponding corner-sum matrix by one.
    
    Another notable object is the \textit{monotone triangle}, which is a triangular array where row $1,...,n$ has $1,...,n$ entries respectively, the numbers in the bottom row are 1 to n (in order), the numbers in each row are strictly increasing from left to right, and the numbers along diagonals are weakly increasing from left to right (see \cite{propp}). Again, by adding (or subtracting) the sub-matrix in \textit{theorem 1} to an alternating sign matrix, we are decreasing (or increasing) one entry in the corresponding monotone triangle by one.
    
    Now that we have proven \textit{theorem 2} for alternating sign triangle, it is possible that there are objects analogous to corner-sum matrix and monotone triangle for alternating sign triangle. Additionally, these objects may also allow us to give a complete characterization of alternating sign triangles.
    
	\justify


\begin{thebibliography}{9}
	
		\bibitem{robbins} 
		D. Robbins and H. Rumsey. “Determinants and Alternating-Sign Matrices,”
		\textit{Advances in Mathematics}. \textbf{62} (1986), 169-184.
		
		\bibitem{brualdi}
		R.A. Brualdi, K.P. Kiernan, S.A. Meyer, M.W. Schroeder. "Patterns of alternating sign matrices," \textit{Linear Algebra Appl}. \textbf{438} (10) (2013) 3967–3990.
		
		\bibitem{propp} 
		J. Propp. “The Many Faces of Alternating-Sign Matrices,”
		\textit{Discrete Mathematics and Theoretical Computer Science Proceedings}.
		\textbf{AA (DM-CCG)}
		(2001), 43-58.
		
		\bibitem{carroll} 
		G. Carroll and D. Speyer. “The Cube Recurrence,”
		\href{https://arxiv.org/abs/math/0403417}{arXiv:math/0403417}
		
		\bibitem{son}
		S. Nguyen. "Characterizing Alternating Sign Triangles," \href{https://arxiv.org/abs/2105.03969}{arXiv:2105.03969}
		
	\end{thebibliography}
\end{document}